\documentclass[10pt]{amsart}
\usepackage{amscd}
\usepackage{amssymb}

\newtheorem{thm}[equation]{Theorem}
\newtheorem{cor}[equation]{Corollary}

\theoremstyle{definition}

\newtheorem{rem}[equation]{Remark}
\newtheorem{exa}[equation]{Example}

\newcommand{\ar}{\rightarrow}
\newcommand{\opn}{\operatorname}
\newcommand{\cat}[1]{\operatorname{\mathsf{#1}}}

\newcommand{\msf}[1]{\mathsf{#1}}
\newcommand{\mbf}[1]{\mathbf{#1}}
\newcommand{\mrm}[1]{\mathrm{#1}}
\newcommand{\mbb}[1]{\mathbb{#1}}

\newcommand{\bsym}[1]{\boldsymbol{#1}}

\begin{document}
\title[Derived Picard Group]
{The Derived Picard Group is a Locally Algebraic Group}
\author{Amnon Yekutieli}
\address{Department of Mathematics,
Ben Gurion University, Be'er Sheva 84105, ISRAEL}
\email{amyekut@math.bgu.ac.il \newline \indent 
http://www.math.bgu.ac.il/$\sim$amyekut}
\date{23.12.00}
\thanks{Partially supported by a grant from the US-Israel Binational
Science Foundation} 

\begin{abstract}
Let $A$ be a finite dimensional algebra over an algebraically 
closed field $\mbb{K}$. The derived Picard group 
$\opn{DPic}_{\mbb{K}}(A)$ is the group of two-sided tilting 
complexes over $A$ modulo isomorphism. 
We prove that $\opn{DPic}_{\mbb{K}}(A)$ is a 
locally algebraic group, and its identity component is 
$\opn{Out}^{0}_{\mbb{K}}(A)$. If $B$ is a derived Morita 
equivalent algebra then 
$\opn{DPic}_{\mbb{K}}(A) \cong \opn{DPic}_{\mbb{K}}(B)$
as locally algebraic groups. Our results extend, and are based on, 
work of Huisgen-Zimmermann, Saor{\'\i}n and Rouquier.
\end{abstract}

\maketitle


Let $A$ and $B$ be associative algebras with $1$ over a field 
$\mbb{K}$. We denote by $\msf{D}^{\mrm{b}}(\cat{Mod} A)$ the bounded 
derived category of left $A$-modules. Let $B^{\circ}$ be the opposite 
algebra, so an $A \otimes_{\mbb{K}} B^{\circ}$-module is a 
$\mbb{K}$-central $A$-$B$-bimodule.
A {\em two-sided tilting complex} over $(A, B)$ is a complex 
$T \in \msf{D}^{\mrm{b}}(\cat{Mod} A \otimes_{\mbb{K}} B^{\circ})$
such there exists a complex 
$T^{\vee} \in 
\msf{D}^{\mrm{b}}(\cat{Mod} B \otimes_{\mbb{K}} A^{\circ})$
and isomorphisms of the derived tensor products
$T \otimes^{\mrm{L}}_{B} T^{\vee} \cong A$
and
$T^{\vee} \otimes^{\mrm{L}}_{A} T \cong B$.
Two-sided tilting complexes were introduced by Rickard in \cite{Rd}.

When $B = A$ we write
$A^{\mrm{e}} := A \otimes_{\mbb{K}} A^{\circ}$.
The set
\[ \opn{DPic}_{\mbb{K}}(A) := 
\frac{ \{ \text{two-sided tilting complexes }
T \in \msf{D}^{\mrm{b}}(\cat{Mod} A^{\mrm{e}}) \}}
{\text{isomorphism}} \]
is the {\em derived Picard group of $A$} (relative to $\mbb{K}$). 
The identity element is the class of $A$, the multiplication is 
$(T_{1}, T_{2}) \mapsto T_{1} \otimes^{\mrm{L}}_{A} T_{2}$, 
and the inverse is 
$T \mapsto T^{\vee} = \opn{RHom}_{A}(T, A)$. 

Denote by $\opn{Out}_{\mbb{K}}(A)$ the group of outer 
$\mbb{K}$-algebra automorphism of $A$, and by 
$\opn{Pic}_{\mbb{K}}(A)$ the Picard group of $A$ (the group of 
invertible bimodules modulo isomorphism). Then there are inclusions
\[ \opn{Out}_{\mbb{K}}(A) \subset \opn{Pic}_{\mbb{K}}(A)
\subset \opn{DPic}_{\mbb{K}}(A) . \]
The first inclusion sends the automorphism $\sigma$ to the 
invertible bimodule $A^{\sigma}$ where the right action is twisted 
by $\sigma$. The second inclusion corresponds to the full embedding
$\cat{Mod} A^{\mrm{e}} \subset 
\msf{D}^{\mrm{b}}(\cat{Mod} A^{\mrm{e}})$.
See \cite{Ye} for details. 

To simplify notation we use the same symbol to 
denote an automorphism $\sigma \in \opn{Aut}_{\mbb{K}}(A)$ and its 
class in $\opn{Out}_{\mbb{K}}(A)$. Likewise for a two-sided 
tilting complex $T$ and its class in $\opn{DPic}_{\mbb{K}}(A)$. 
The precise meaning is always clear from the context. 

Now assume $\mbb{K}$ is algebraically closed and $A$ is a 
finite dimensional $\mbb{K}$-algebra. Then the group 
$\opn{Aut}_{\mbb{K}}(A) = \opn{Aut}_{\cat{Alg} \mbb{K}}(A)$ 
of $\mbb{K}$-algebra automorphisms is a 
linear algebraic group, being a closed subgroup of 
$\opn{GL}(A) = \opn{Aut}_{\cat{Mod} \mbb{K}}(A)$. 
This induces a structure of linear algebraic group on the quotient
$\opn{Out}_{\mbb{K}}(A)$. Denote by $\opn{Out}^{0}_{\mbb{K}}(A)$
the identity component. 

Examples calculated in \cite{MY} indicated that the whole group 
$\opn{DPic}_{\mbb{K}}(A)$ should carry a geometric structure (cf.\ 
Example \ref{exa1} below). This is our first main result Theorem 
\ref{thm1}. 

A result of Brauer says that the group 
$\opn{Out}^{0}_{\mbb{K}}(A)$ is a Morita invariant of $A$: if $A$ 
and $B$ are Morita equivalent $\mbb{K}$-algebras then 
$\opn{Out}^{0}_{\mbb{K}}(A) \cong \opn{Out}^{0}_{\mbb{K}}(B)$. 
In \cite{HS} and \cite{Ro} this is extended to derived 
Morita equivalence. Our Theorem \ref{thm2} extends these 
results further. 

We shall need the following variant of the result of 
Huisgen-Zimmermann, Saor{\'\i}n and Rouquier.

\begin{thm} \label{thm3}
Let $A$ and $B$ be finite dimensional $\mbb{K}$-algebras. Suppose
$T \in$ \linebreak
$\msf{D}^{\mrm{b}}(\cat{Mod} A \otimes_{\mbb{K}} B^{\circ})$
is a two-sided tilting complex over $(A, B)$, with inverse 
$T^{\vee} \in 
\msf{D}^{\mrm{b}}(\cat{Mod} B \otimes_{\mbb{K}} A^{\circ})$. 
Then for any element
$\sigma \in \opn{Out}^{0}_{\mbb{K}}(A)$
the two-sided tilting complex
\[ \phi_{T}^{0}(\sigma) := 
T^{\vee} \otimes^{\mrm{L}}_{A} A^{\sigma} \otimes^{\mrm{L}}_{A} T
\in \opn{DPic}_{\mbb{K}}(B) \]
is in $\opn{Out}^{0}_{\mbb{K}}(B)$.
The group homomorphism
\[ \phi_{T}^{0}: \opn{Out}_{\mbb{K}}^{0}(A) \ar 
\opn{Out}_{\mbb{K}}^{0}(B) \]
is an isomorphism of algebraic groups.
\end{thm}

\begin{proof}
According to \cite[Theorem 17]{HS} or \cite[Th\'eor\`eme 4.2]{Ro} 
there is an isomorphism of algebraic groups
$\phi^{0}: \opn{Out}^{0}_{\mbb{K}}(A) \ar 
\opn{Out}^{0}_{\mbb{K}}(B)$
induced by $T$. Letting 
$\tau := \phi^{0}(\sigma) \in \opn{Out}^{0}_{\mbb{K}}(B)$
one has 
\[ T \otimes_{B} B^{\tau} \cong
A^{\sigma} \otimes_{A} T 
\text{ in } \msf{D}(\cat{Mod} A \otimes_{\mbb{K}} B^{\circ})
. \]
Applying $T^{\vee} \otimes^{\mrm{L}}_{A} -$ to this isomorphism
we see that $B^{\tau} \cong \phi^{0}_{T}(\sigma)$
in $\msf{D}(\cat{Mod} B^{\mrm{e}})$, so 
$\tau = \phi^{0}_{T}(\sigma)$ in $\opn{DPic}_{\mbb{K}}(B)$. We 
conclude that $\phi^{0}_{T} = \phi^{0}$.
\end{proof}

A {\em locally algebraic group} over $\mbb{K}$ is a group $G$, 
with a normal subgroup $G^{0}$, such that $G^{0}$ is a connected 
algebraic group over $\mbb{K}$, each coset of $G^{0}$ is a 
variety, and multiplication and inversion are morphisms of 
varieties. 
A morphism $\phi: G \to H$ of locally algebraic groups is a group 
homomorphism such that $\phi(G^{0}) \subset H^{0}$ and the 
restriction $\phi^{0}: G^{0} \to H^{0}$ is a morphism of varieties. 
We call $\phi$ an open immersion if $\phi$ is injective and
$\phi^{0}$ is an isomorphism.

In other words $G$ is the group of rational points 
$\bsym{G}(\mbb{K})$ of a reduced group scheme $\bsym{G}$ locally of 
finite type over $\mbb{K}$, in the sense of \cite[Expos\'e 
VI$_{\mrm{A}}$]{SGA3}. A morphism $\phi: G \to H$ corresponds to a 
morphism $\bsym{\phi}: \bsym{G} \to \bsym{H}$ of group schemes over 
$\mbb{K}$. 

Here is our first main result.

\begin{thm} \label{thm1}
Let $A$ be a finite dimensional $\mbb{K}$-algebra. 
Then the derived Picard group 
$\opn{DPic}_{\mbb{K}}(A)$ is a locally algebraic group over
$\mbb{K}$. The inclusion 
$\opn{Out}_{\mbb{K}}(A) \subset \opn{DPic}_{\mbb{K}}(A)$ 
is an open immersion.
\end{thm}

In particular the identity components coincide:
$\opn{Out}^{0}_{\mbb{K}}(A) = \opn{DPic}^{0}_{\mbb{K}}(A)$.

\begin{proof}
Theorem \ref{thm3} with $A = B$ implies that the subgroup
$\opn{Out}^{0}_{\mbb{K}}(A) \subset \opn{DPic}_{\mbb{K}}(A)$
is normal, and for any two-sided tilting complex $T$
the conjugation 
$\phi^{0}_{T} : \opn{Out}^{0}_{\mbb{K}}(A) \to
\opn{Out}^{0}_{\mbb{K}}(A)$
is an automorphism of algebraic groups.  

Let us now switch to the notation $T_{1} \cdot T_{2}$ and $T^{-1}$ 
for the operations in $\opn{DPic}_{\mbb{K}}(A)$. 
Define an algebraic variety structure on each coset 
$C = T \cdot \opn{Out}^{0}_{\mbb{K}}(A) \subset 
\opn{DPic}_{\mbb{K}}(A)$
using the multiplication map
$P \mapsto T \cdot P$, 
$P \in \opn{Out}^{0}_{\mbb{K}}(A)$. 
Since $\phi^{0}_{T}$ is an automorphism of algebraic groups, the 
variety structure is independent of the representative $T \in C$. 

Let us prove that $\opn{DPic}_{\mbb{K}}(A)$
is a locally algebraic group. For 
$P_{1}, P_{2} \in \opn{Out}^{0}_{\mbb{K}}(A)$ 
and
$T_{1}, T_{2} \in \opn{DPic}_{\mbb{K}}(A)$,
multiplication is the morphism
\[ (T_{1} \cdot P_{1}) \cdot (T_{2} \cdot P_{2}) =
(T_{1} \cdot T_{2}) \cdot (\phi^{0}_{T_{2}}(P_{1}) \cdot P_{2})
. \]
Similarly for the inverse:
\[ (T \cdot P)^{-1} = T^{-1} \cdot \phi^{0}_{T}(P)^{-1} . \]
\end{proof}

\begin{exa} \label{exa1}
Let $\vec{\Omega}_{n}$ be the quiver with two vertices $x, y$ and 
$n$ arrows $x \xrightarrow{\alpha_{i}} y$. Let $A$ be 
the path algebra $\mbb{K} \vec{\Omega}_{n}$. According to 
\cite[Theorem 5.3]{MY}, 
$\opn{Out}_{\mbb{K}}(A) \cong \opn{Pic}_{\mbb{K}}(A) 
\cong \opn{PGL}_{n}(\mbb{K})$
and
\[ \opn{DPic}_{\mbb{K}}(A) \cong \mbb{Z} \times 
\bigl( \mbb{Z} \ltimes \opn{PGL}_{n}(\mbb{K}) \bigr) . \]
In the semi-direct product a generator $T$ of $\mbb{Z}$ acts on 
a matrix $\sigma \in \opn{PGL}_{n}(\mbb{K})$ by 
$\phi^{0}_{T}(\sigma) = (\sigma^{-1})^{\mrm{t}}$.
This is clearly a morphism of varieties, so 
$\opn{DPic}_{\mbb{K}}(A)$ is indeed a locally algebraic group.
\end{exa}

Our second main result relates two algebras. Recall that the 
algebras $A$ and $B$ are derived Morita equivalent over $\mbb{K}$ 
if there is a $\mbb{K}$-linear equivalence of triangulated 
categories
$\msf{D}^{\mrm{b}}(\cat{Mod} A) \approx 
\msf{D}^{\mrm{b}}(\cat{Mod} B)$. 

\begin{thm} \label{thm2}
Suppose $A$ and $B$ are two finite dimensional $\mbb{K}$-algebras, 
and assume they are derived Morita equivalent over $\mbb{K}$. Then 
$\opn{DPic}_{\mbb{K}}(A) \cong \opn{DPic}_{\mbb{K}}(B)$
as locally algebraic groups.
\end{thm}

\begin{proof}
It is known that there exist two-sided tilting complexes
$T \in \msf{D}(\cat{Mod} A \otimes_{\mbb{K}} B^{\circ})$; choose 
one. We obtain a group isomorphism
\[ 
\phi_{T}: \begin{cases}
\opn{DPic}_{\mbb{K}}(A) \to \opn{DPic}_{\mbb{K}}(B) , \\
S \mapsto 
T^{\vee} \otimes^{\mrm{L}}_{A} S \otimes^{\mrm{L}}_{A} T .
\end{cases} \]
By Theorem \ref{thm3}, $\phi_{T}$ restricts to an isomorphism of 
algebraic groups
$ \phi_{T}^{0} : \opn{Out}^{0}_{\mbb{K}}(A) \to 
\opn{Out}^{0}_{\mbb{K}}(B)$.
So $\phi_{T}$ is an isomorphism of locally algebraic groups.
\end{proof}

We end the paper with a corollary and some remarks.
Suppose $\msf{C}$ is a $\mbb{K}$-linear triangulated 
category that's equivalent to a small category. Denote by 
$\opn{Out}^{\mrm{tr}}_{\mbb{K}}(\msf{C})$
the group of $\mbb{K}$-linear triangle auto-equivalences of 
$\msf{C}$ modulo natural isomorphism. Let $\cat{mod} A$ stand for
the category of finitely generated $A$-modules. 

\begin{cor} \label{cor1}
Suppose $\msf{C}$ is a $\mbb{K}$-linear triangulated category that 
is equivalent to $\msf{D}^{\mrm{b}}(\cat{mod} A)$ for some 
hereditary finite dimensional $\mbb{K}$-algebra $A$. Then 
$\opn{Out}^{\mrm{tr}}_{\mbb{K}}(\msf{C})$ is a locally algebraic 
group. 
\end{cor}

\begin{proof}
Trivially 
$\opn{Out}^{\mrm{tr}}_{\mbb{K}}(\msf{C}) \cong
\opn{Out}^{\mrm{tr}}_{\mbb{K}}(\msf{D}^{\mrm{b}}(\cat{mod} A))$, 
and by \cite[Corollary 0.11]{MY} we have
$\opn{Out}^{\mrm{tr}}_{\mbb{K}}(\msf{D}^{\mrm{b}}(\cat{mod} A))
\cong \opn{DPic}_{\mbb{K}}(A)$.
\end{proof}

\begin{exa} \label{exa2}
Beilinson \cite{Be} proved that 
$\msf{D}^{\mrm{b}}(\cat{Coh} \mbf{P}^{1}_{\mbb{K}}) \approx
\msf{D}^{\mrm{b}}(\cat{mod} \mbb{K} \vec{\Omega}_{2})$,
where \linebreak
$\cat{Coh} \mbf{P}^{1}_{\mbb{K}}$ is the category of 
coherent sheaves on the projective line, and $\vec{\Omega}_{2}$ is 
the quiver from Example \ref{exa1}. Therefore 
$\opn{Out}^{\mrm{tr}}_{\mbb{K}}
(\msf{D}^{\mrm{b}}(\cat{Coh} \mbf{P}^{1}_{\mbb{K}}))$
is a locally algebraic group. This should be compared to  
Remark \ref{rem1} below; see also \cite[Remark 5.4]{MY}.
\end{exa}

\begin{rem} \label{rem1}
Suppose $X$ is a smooth projective variety over $\mbb{K}$ with 
ample canonical or anti-canonical bundle. Bondal and Orlov \cite{BO}
prove that 
\[ \opn{Out}^{\mrm{tr}}_{\mbb{K}}(\msf{D}^{\mrm{b}}(\cat{Coh} X))
\cong \left( \opn{Aut}_{\mbb{K}}(X) \ltimes \opn{Pic}(X) \right)
\times \mbb{Z} . \]
Here $\opn{Pic}(X)$ is the group of line bundles. Thus 
$\opn{Out}^{\mrm{tr}}_{\mbb{K}}(\msf{D}^{\mrm{b}}(\cat{Coh} X))
\cong G \times D$
where $G$ is an algebraic group and $D$ is a discrete group, and 
in particular this is a locally algebraic group.
\end{rem}

\begin{rem} \label{rem2}
In \cite{Or}, Orlov gives an example of an abelian variety over 
$\mbb{K}$ such that 
\[ \opn{Out}^{\mrm{tr}}_{\mbb{K}}(\msf{D}^{\mrm{b}}(\cat{Coh} X))
\cong D \ltimes (X \times \widehat{X})(\mbb{K}) , \]
where $D$ is a discrete group (an extension of 
$\opn{SL}_{2}(\mbb{Z})$ by $\mbb{Z}$) and $\widehat{X}$ is the dual 
abelian variety. The group $D$ acts (nontrivially) via
$\opn{Aut}_{\mbb{K}}(X  \times \widehat{X})$
and hence 
$\opn{Out}^{\mrm{tr}}_{\mbb{K}}(\msf{D}^{\mrm{b}}(\cat{Coh} X))$
is a locally algebraic group.
\end{rem}

\medskip \noindent \textbf{Acknowledgments.}\
I wish to thank Birge Huisgen-Zimmermann for showing me \cite{HS}, 
and for very helpful comments on an earlier version of this paper. 
Thanks also to Raphael Rouquier for calling my attention to \cite{Ro} 
and for illuminating discussions.

\end{document}